\documentclass{article}
\usepackage{amsmath,amsthm}
\usepackage{graphicx,epstopdf,epsfig,multirow,epic,bm}

\normalsize \rm
\begin{document}

\begin{center}
{\Large  \textbf {An extremal problem: How small scale-free graph can be  }}\\[12pt]
{\large Fei Ma$^{a,}$\footnote{~E-mail: mafei123987@163.com. },\quad Ping Wang$^{b,c,d,}$\footnote{~E-mail: pwang@pku.edu.cn.} \quad  and  \quad  Bing Yao$^{a,}$\footnote{~E-mail: yybb918@163.com.} }\\[6pt]
{\footnotesize $^{e}$ School of Electronics Engineering and Computer Science, Peking University, Beijing 100871, China\\
$^{b}$ School of Software and Microelectronics, Peking University, Beijing  102600, China\\
$^{c}$ National Engineering Research Center for Software Engineering, Peking University, Beijing, China\\
$^{d}$ Key Laboratory of High Confidence Software Technologies (PKU), Ministry of Education, Beijing, China\\
$^{e}$ College of Mathematics and Statistics, Northwest Normal University, 730070 Lanzhou, China.}\\[12pt]
\end{center}

\begin{quote}
\textbf{Abstract:} The bloom of complex network study, in particular, with
respect to scale-free ones, is considerably triggering the research
of scale-free graph itself. Therefore, a great number
of interesting results have been reported in the past, including
bounds of diameter. In this paper, we focus mainly on a
problem of how to analytically estimate the lower bound of
diameter of scale-free graph, i.e., how small scale-free graph can be. Unlike some pre-existing
methods for determining the lower bound of diameter, we
make use of a constructive manner in which one candidate
model $\mathcal{G^{*}}(\mathcal{V^{*}},\mathcal{E^{*}})$ with ultra-small diameter can be generated.
In addition, with a rigorous proof, we certainly demonstrate
that the diameter of graph $\mathcal{G^{*}}(\mathcal{V^{*}},\mathcal{E^{*}})$ must be the
smallest in comparison with that of any scale-free graph.
This should be regarded as the tight lower bound.


\end{quote}

\vskip 1cm

\section{Introduction}

As a special member of graph family, scale-free graph has attracted considerable attention in the past. One of important reasons for this is the bloom of complex network study in the last two decades, in particular, in terms of scale-free ones \cite{Albert-1999}. In the jargon of graph theory, we let $\mathcal{G}(\mathcal{V},\mathcal{E})$ denote a graph where $\mathcal{V}$ and $\mathcal{E}$ represent, respectively, vertex set and edge set. Accordingly, the symbols $|\mathcal{V}|$ and $|\mathcal{E}|$ are the order and size of graph $\mathcal{G}(\mathcal{V},\mathcal{E})$, separately. Mathematically, diameter of a graph $\mathcal{G}(\mathcal{V},\mathcal{E})$, denoted by $D$, is the maximum over distances of all possible vertex pairs. For a pair of vertices $u$ and $v$, distance between them, denoted by $d_{uv}$, is the edge number of any shorted path joining vertex $u$ and $v$. Here we just consider simple graphs, i.e., one type of graphs that have no multi-edges and loops.

\subsection{Description of problems}

Given a graph $\mathcal{G}(\mathcal{V},\mathcal{E})$, one can determine with degree distribution, usually called degree sequence, whether it is scale-free or not. Graph $\mathcal{G}(\mathcal{V},\mathcal{E})$ can be considered scale-free if its degree distribution $P(k)$ follows

\begin{equation}\label{MF-1}
P(k)\sim k^{-\gamma}, \quad 1<\gamma
\end{equation}
where $P(k)$ is the probability of randomly selecting a vertex with degree equal to $k$ from graph $\mathcal{G}(\mathcal{V},\mathcal{E})$. In discrete case, Eq.(\ref{MF-1}) can be expressed in an alterative manner as follows

\begin{equation}\label{MF-2}
P_{cum}(k_{i}\geq k)=\frac{\sum_{k_{i}\geq k}N_{k_{i}}}{|\mathcal{V}|}\sim k^{-(\gamma+1)}
\end{equation}
in which $N_{k_{i}}$ is the total number of degree $k_{i}$ vertices. It is convention to call $P_{cum}(k_{i}\geq k)$ the accumulative degree distribution of graph $\mathcal{G}(\mathcal{V},\mathcal{E})$. It is worth noting that we will prove the main result of this paper, i.e., how small scale-free graphs can be \footnote{Commonly, using diameter estimates whether a scale-free graph is small or not. Here, we also employ such a topological parameter to quantify graphs with scale-free feature.},
 using Eq.(\ref{MF-2}) mainly because our graph is constructed in a deterministic manner.

\subsection{Related work}

Before starting with our discussions, we need to recall some previous related work in this field.

\textbf{Result 1} In \cite{Cohen-2003}, the authors demonstrated using analytical arguments that scale-free graphs with $2<\gamma<3$ have a much smaller diameter, behaving as $D\sim\ln\ln N$. For $\gamma=3$, this yields $D\sim \ln N/\ln\ln N$, while for $\gamma>3$, $D \sim\ln N$.

\textbf{Result 2} In \cite{Bollobas-2004}, the authors showed that fixing an integer $m\geq2$ and a positive real number $\epsilon$, then a.e. $G^{n}_{m}\in \mathcal{G}^{n}_{m}$ is connected and has diameter $D(G^{n}_{m})$ satisfying

$$(1-\epsilon)\log n/\log\log n\leq D(G^{n}_{m}) \leq (1+\epsilon)\log n/\log\log n.$$
where $n$ is the order of graph $G^{n}_{m}$.

\textbf{Result 3} In \cite{Chung-2002}, the authors stated that for a random sparse graph $\mathcal{G}(\mathcal{V},\mathcal{E})$ with admissible expected degree sequence $(w_{1},\dots,w_{|\mathcal{V}|})$, the diameter is almost surely

$$\Theta(\log |\mathcal{V}|/\log \widetilde{d})$$
here $\widetilde{d}$ is second order average degree of graph $\mathcal{G}(\mathcal{V},\mathcal{E})$ that satisfies $0< \log \widetilde{d} \ll \log |\mathcal{V}|$.

As shown above, for all scale-free graphs, including deterministic and stochastic, or sparse and dense, the tight lower bound of their diameters seems to not be obtained. To do this, in this paper, we will generate a deterministic graph $\mathcal{G^{*}}(\mathcal{V^{*}},\mathcal{E^{*}})$ in a concise fashion. And then, it turns out to both be scale-free and have smallest diameter. First, below is the theorem whose complete proof will be deferred in Section 3.

\textbf{Theorem } For all scale-free graphs, the tight lower bound of the diameter is able to be equal to $2$.

From now then, let us turn our insight into construction of candidate graph $\mathcal{G^{*}}(\mathcal{V^{*}},\mathcal{E^{*}})$.

\section{Construction}

\begin{figure}
\centering
  \includegraphics[height=5cm]{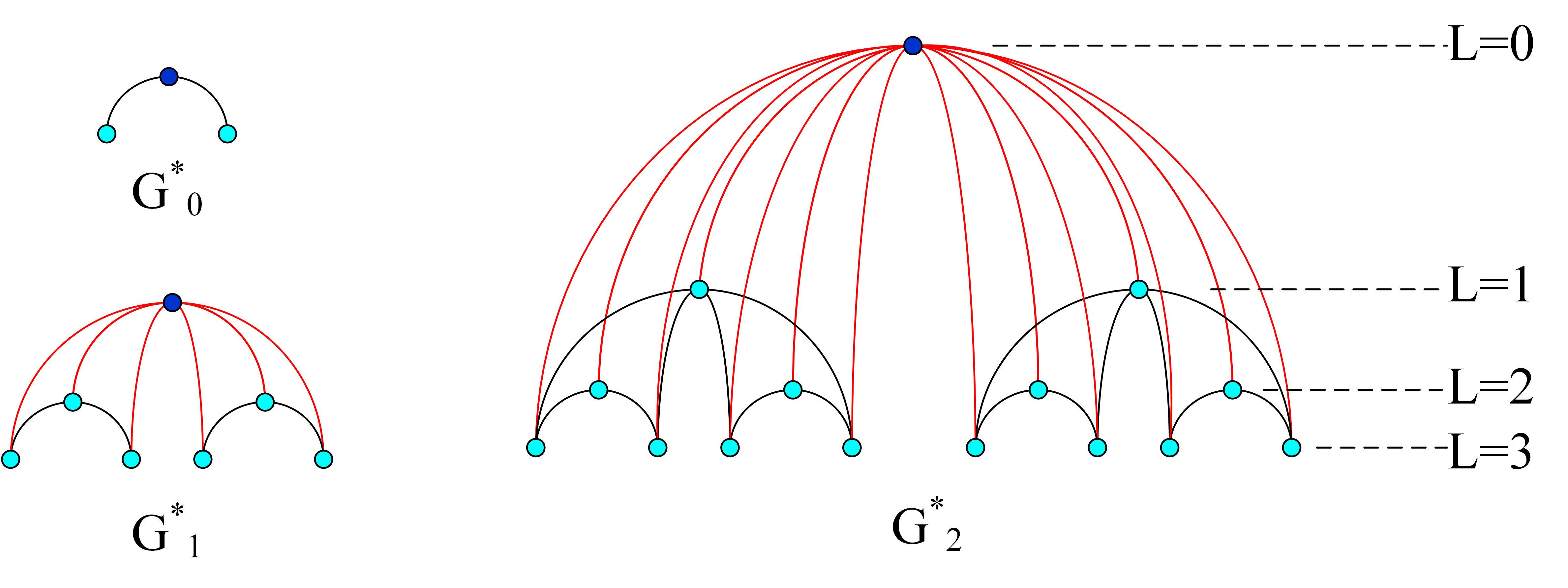}\\
{\small Fig.1. The diagram of first three examples of graphs $\mathcal{G^{*}}_{t}$.  }
\end{figure}

This section will introduce our proposed graph $\mathcal{G^{*}}(\mathcal{V^{*}},\mathcal{E^{*}})$ and provide a rigorous proof of Theorem.

First, the seed, denoted by $\mathcal{G^{*}}_{0}$, is a star with $2$ leaves as shown in the top-left panel in Fig.1.
The second graph $\mathcal{G^{*}}_{1}$ can be generated based on $\mathcal{G^{*}}_{0}$ in the following manner,

(\textbf{i}) making $2$ duplications of seed $\mathcal{G^{*}}_{0}$ labelled as $\mathcal{G^{*}}_{0}^{i}$,

(\textbf{ii}) taking an active vertex (blue online),

(\textbf{iii}) connecting that active vertex to each vertex in star $\mathcal{G^{*}}_{0}^{i}$.

For our purpose, we group all vertices of $\mathcal{G^{*}}_{1}$ into three classes, i.e., that active vertex allocated at the level $0$, denoted by $L=0$, the centers of stars $\mathcal{G^{*}}_{0}^{i}$ at the level $L=1$ and the remaining vertices of graph $\mathcal{G^{*}}_{1}$ at the level $L=2$. With the help of such a classification, the next graph $\mathcal{G^{*}}_{2}$ can be generated from $\mathcal{G^{*}}_{1}$ using both the above methods \textbf{i-ii} and an additional technique as follows

(\textbf{iv}) connecting that active vertex to each vertex in star $\mathcal{G^{*}}_{1}^{i}$ and simultaneously deleting all edges but for those adjacent to vertices at the level $L=2$ in $\mathcal{G^{*}}_{1}^{i}$ .

Since then,  for time step $t\geq3$, the young graph $\mathcal{G^{*}}_{t}$ can be built on the basis of $2$ duplications of the preceding graph $\mathcal{G^{*}}_{t-1}$ using procedures \textbf{i-ii} and \textbf{iv}. To be more concrete, the graph $\mathcal{G^{*}}_{2}$ is illustrated in the rightmost panel of Fig.1. As will be shown shortly, graph $\mathcal{G^{*}}_{t}$ is a candidate model.

\section{Proof of theorem}

In view of the construction of graph, it is easy to calculate the order $\mathcal{V^{*}}_{t}$ and size $\mathcal{E^{*}}_{t}$ of $\mathcal{G^{*}}_{t}$ in the following form

$$|\mathcal{V^{*}}_{t}|=2^{t+2}-1, \quad |\mathcal{E^{*}}_{t}|=2^{t+1}(t+2)-2.$$

Similarly, we can obtain a list consisting of degree sequence of graph $\mathcal{G^{*}}_{t}$

\begin{center}
\begin{tabular}{c|c|c|c|c|c|c|c|cccc}
  \hline
  $k_{t_{i},t}$ & $2^{t+2}$ & $2^{t}+1$ & $...$ & $2^{t_{i}}+1$ & $...$ & $2^{2}+1$ & $2+1$ & $t+1$ \\
   \hline
 $N_{t_{i},t}$ & $1$ & $2$ & $...$ & $2^{t-t_{i}+1}$ & $...$ & $2^{t-1}$ & $2^{t}$& $2^{t+1}$ \\
  \hline
\end{tabular}
\end{center}

Using Eq.(\ref{MF-2}), we can have

\begin{equation}\label{MF-3}
P_{cum}(k\geq k_{t_{i},t})=\left\{
\begin{array}{ll}
\quad k_{t_{i},t}^{-\gamma_{\alpha}}\; , \qquad & \quad k_{t_{i},t}> t+1\\
k_{t_{i},t}^{-\gamma_{\alpha}}+\frac{1}{2}\; , & \quad k_{t_{i},t}\leq t+1 \\
\end{array}
\right.
\end{equation}
Taking the derivative of both sides in Eq.(\ref{MF-3}) with respect to $k$ produces
\begin{equation}\label{MF-4}
P(k)\sim k^{-\gamma}, \quad \gamma=\gamma_{\alpha}+1=2.
\end{equation}
This implies that graph $\mathcal{G^{*}}_{t}$ follows power-law distribution.

Now, let us prove Theorem in the end of Section 1.

\emph{Proof } Consider a connected scale-free graph $\mathcal{G}(\mathcal{V},\mathcal{E})$ with diameter $D$, if $D$ might be equivalent to $1$ then it means that for arbitrary vertex pairs $u$ and $v$ there must be $1\leq d_{uv}\leq D$. This surely suggests that graph $\mathcal{G}(\mathcal{V},\mathcal{E})$ is a complete graph and can not be scale-free. So clear to see that $D$ is strictly larger than $1$. With our graph $\mathcal{G^{*}}_{t}$ with scale-free feature, it can easily determine that the diameter $D^{*}$ is exactly equal to $2$. This completes our proof.

In a word, the lower bound of diameter of scale-free graph may be an ultra-small constant $2$.

\end{document}